\documentclass[12pt]{article}
 \usepackage{amsmath,amsfonts,theorem}
 \numberwithin{equation}{section}

 \newtheorem{prop}{Proposition}[section]

 \newtheorem{cor}{Corollary}[section]

 \theorembodyfont{\rmfamily}

 \newcommand{\qed}{\ifhmode\unskip\nobreak\fi\quad\ensuremath\square}

 \newcommand{\sL}{\mathcal L}

 \newcommand{\Oh}{\mathcal O}

 \newcommand{\al}{\alpha}
 \newcommand{\be}{\beta}

 \newcommand{\ga}{\gamma}
 
 \newcommand{\om}{\omega}
 
 \newcommand{\De}{\Delta}

 \newcommand{\la}{\lambda}

 \newcommand{\PP}{\mathbb P}
 \newcommand{\C}{\mathbb C}

 \newcommand{\R}{\mathbb R}

 \renewcommand{\Im}{\operatorname{Im}}

 \newcommand{\Sing}{\operatorname{Sing}}
 
\begin{document}

 \title{Special lagrangian fibrations on flag variety $F^3$}

 \author{Nikolay A. Tyurin\footnote{BLTPh JINR (Dubna) and MGU PS (MIIT).  This work
 was partially supported by RFBR, grants
N 08-01-00095, N 08-01-00392a}}
\date{}
 \maketitle

\begin{center}
 {\em                                  }
 \end{center}
 \bigskip

\section*{Introduction}

Lagrangian fibrations of symplectic manifolds are important for several
applications.

From the classical mechanics point of view if a phase space is fibered on
lagrangian submanifolds then the motion of the corresponding system with a
Hamiltonian function whose hamiltonian vector field is parallel to the
fibration can be described in terms of the fibers, and this can be done even in
the case when certain fibers are singular.

An approach adopted in Geometric Quantization is realized in the same setup ---
if a given phase space is fibered on lagrangian submanifolds then one can take
the so - called Bohr - Sommerfeld fibers which form a discrete subset, and this
subset generates the corresponding Hilbert space. Here one should mention that
the notion of Bohr - Sommerfeld lagrangian submanifold is valid for certain
singular lagrangian fibers.

Homological Mirror Symmetry today can be treated by considerations of  a
lagrangian fibration and derivation a finite set of fibers which satisfy
certain condition in terms of Fukaya - Floer cohomology. Again, one could try
to count holomorphic discs for certain singular tori, so it seems that this
approach to HMS can be extended to the cases when the fibrations have singular
fibers.

The ideal and model situation for the story originates in toric geometry,
\cite{1}. If one works with a symplectic toric manifold  which is the phase
space of a completely integrable mechanical system then there is a canonical
lagrangian fibration of the manifold parameterized by a convex polytop and this
leads to the introduction of the action - angle variables with complete
solution of the system. At the same time in the toric case the methods for GQ
and HMS mentioned above work succesfully. But the problem is that the set of
toric manifolds is not too big and there are many other manifolds which do not
admit a toric structure due to topological reasons but which must be involved
in the consideration.

For a compact symplectic toric manifold $(M, \om)$ the canonical lagrangian
fibration is not completely regular --- a dimensional reduction takes place so
smooth lagrangian tori fill not whole $M$ but the complement $M \ D$ where $D$
is a "symplectic divisor" whose homology class is Poincare dual to the
anticanonical class of $M$. This symplectic submanifold $D$ of real codimension
2 consists of a number of irreducible smooth components $D_i$, and each
component is as a symplectic manifold fibered on lagrangian tori which are
isotorpical being considered in $M$. Totally it can be described in terms of
the moment polytop $P$: each facet of it is a moment polytop of a smooth
symplectic toric manifold of real dimension equals to twice the dimension of
the facet. In this setup it is reasonable to call $D$ {\it the boundary
divisor}.

Let us extend the class of symplectic toric manifolds. Here and in other
papers, \cite{3}, \cite{4}, we study compact smooth symplectic manifolds which
admit lagrangian fibrations of almost the same type as in the toric case unless
the allowance of the existence of singular tori in the fibrations. We require
that generic fiber must be smooth, and the types of the allowed singularties
are exhausted by the following models. If $T^n$ is a smooth $n$ - dimensional
torus, then we can take a decomposition
$$
T^n = T^k \times T^{n-k}
$$
and then the following precedure can be applied to our $T^n$: at a point $p \in
T^{n-k}$ the fiber $T^k_p = T^k \times \{ p \}$ collapses to $k-1$ -
dimensional torus. If we denote the corresponding degeneration operator as
$C^k_p$, then we say that a singular torus which is allowed in our fibrations
must be of the type $(C^{k_1}_{p_1} ... C^{k_l}_{p_l} (T^n)$. Note that this
type of singularities is familiar in classical mechanics --- these correspond
to separatrices; while a smooth torus in the fibration carries loop  type
trajectories, a singular torus can have lines and stable points as solutions.

In the present text we construct lagrangian fibrations of the type described
above for complex flag variety $F^3$, full flag variety for $\C^3$. A
lagrangian fibration of $F^3$ with boundary divisor $D$ formed by 4 del Pezzo
surfaces $\C \PP^2_1$ and with a subfamily of singular tori which has real
dimension one is described below. Generic fiber of the fibration is smooth
lagrangian torus, and the singular fiber is given by $C^2_p$ applied to $T^3$.
This fibration is not unique given by our construction, but we distinguish it
and call it {\it minimal} since it has minimal set of singular tori. The
construction is based on the fact that $F^3$ admits a {\it pseudo toric
structure}, a generalization of the toric structure, proposed in [3] and [4].
Any pseudo toric manifold admits a variety of lagrangian fibrations of the
desired type.

It is well known Gelfand - Zeytlin system on $F^n$, a completely integrable
system, which is not however {\it effective} since $F^n$ is not toric. The
image of $F^n$ under the "action" map is the Gelfand - Zeytlin polytop, and the
difference between the toric case is manifested by the fact that the
dimensional reduction is broken in this case: certain facets of the Gelfand -
Zeytlin polytop are not the images of complex subsets of $F^n$, and they are
not the images of certain submanifolds of right dimensions in $F^n$ (of the
dimensions equal to twice the dimensions of the facets), f.e. there are
vertices whose preimages are not points in $F^n$  (there is a vertex of maximal
degeneration whose preimage is a lagrangian $n$ - dimensional sphere). Thus our
construction in contrast to the Gelfand - Zeytlin picture gives an alternative
way to get a lagrangian fibration on $F^n$ --- we keep the dimensional
reduction argument as in the toric case but we loose the smoothness of the
fibers. Indeed, the boundary divisor for the Gelfand - Zeytlin system on the
flag variety $F^3$ consists of 6 components: 2 complex and 4 non complex, and
there is unique 4 - valent vertex (which lies in the intersection of these 4
components) in the Gelfand - Zeytlin polytop whose preimage is a lagrangian
sphere in $F^3$. In constrast the boundary divisor of our minimal lagrangian
fibration consists of 4 complex components.

Why the minimal lagrangian fibrations of $F^3$  are interesting for us?

First, it gives new special solution to  the classical mechanical problems on
$F^3$. F.e.,  one  notes that three functions $\la^i_j$ of the Gelfand -
Zeytlin system can be combined to give two functions $F_1, F_2$ we are working
with in our construction. This means that if one takes a Hamiltonian $H$ on
$F^3$ of the form
$$
H = \al F_1 + \be F_2
$$
where $\al, \be \in \R$, then the trajectories of the motion generated by $H$,
can be found in two ways: either one takes the fibers of the Gelfand - Zeytlin
lagrangian fibration and gets the motion restricting it to the fibers, or one
can take the minimal fibration and study the fibers with the restricted motion.
The answer, of course, must be the same, but perhaps in certain circumstances
it would be convinient to work with our minimal fibrations, for example, in the
problem of minimal loop searching.

Second, one can apply the standard methods of Geometric Quantization to $F^3$
endowed with the minimal lagrangian fibration. For a singular torus which
appears in the fibration the notion of Bohr - Sommerfeld lagrangian cycle is
still valid; the collapsed torus has a set of "virtual" periods, corresponding
to the shrinked generators of the fundamental group; and these virtual periods
must be trivial with respect to any integer symplectic form. Therefore one can
define the set of Bohr - Sommerfeld fibers (with respect to the anticanonical
class since $F^3$ is a Fano variety); one can prove that the set is discrete
and finite, and thus the corresponding Hilbert space as usual is spanned by
this finite set. The advantage of the minimal lagrangian fiberation is a
realrionship between lagrangian and Kahler quantization for $F^3$ --- as in the
original toric case one has a version of the Danilov - Khovanskiy
correspondince between holomorphic sections of the anticanonical class and the
lagrangian cycle Bohr - Sommerfeld with respect to the anticanonical class.
However for this set of problems we already have a classical result [5] which
solves the Geometric Quantization problem in terms of representation theory.

Third, if one is interested in {\it special lagrangian fibrations} on Fano
varieties, inroduced recently by D. Auroux in \cite{2}, it is clear that if one
studies the flag variety $F^3$ (or more generally $F^n$) then the Gelfand -
Zeytlin fibration is not of the special type, and moreover it doesn't lead to
something close to the desired type. At the same time the minimal lagrangian
fibration on $F^3$  constructed below is special lagrangian, and this is the
point of our main interest.

 The last item is rather specific, but
nevertheless it may lead to certain consistent results in the studies of
lagrangian geometry of Fano varieties.

{\bf Acknowledgements.} I would like to thank S. Belyov,  D. Auroux and D.
Orlov for helpfull discussion and remarks. Special thank goes to Y. Nohara
whose talk at GEOQUANT - 2009 has inspired this work. It follows that I should
thank the organizers of the conference, first of all M. Schlichenmaier and A.
Sergeev, which I do with great pleasure.

\section{Pseudotoric structure on the flag variety $F^3$}

We construct a big family of lagrangian fibrations on $F^3$ using the notion of
pseudotoric structures introduced in \cite{3}, \cite{4}. We have the following

\begin{prop} The flag variety $F^3$ is pseudotoric.
\end{prop}

Recall that pseudotoric structure on a symplectic manifold $(M, \om_M)$ of real
dimension $2n$ is a set of data $(f_1, ..., f_k, B, \psi, (Y, \om_Y))$ where
$f_i$ are smooth Morse functions in involution, $\psi$ defines a family of
symplectic toric submanifolds with the base set $B \subset M$ parameterized by
a symplectic toric manifold $(Y, \om_Y)$  of real dimension $2(n-k)$ so
$$
\psi: M - B \to Y
$$
is a map with generically smooth symplectic fibers, such that the following
compatibility conditions hold:

(1) for each $p \in Y$ the corresponding element of the family
$$
\overline{\psi^{-1}(p)} = \psi^{-1}(p) \cup B
$$
is preserved by the Hamiltonian action of each $f_i$;

(2) for each smooth function $h \in C^{\infty}(Y, \R)$ we have the following
identity
$$
\nabla_{\psi} X_h = \tau X_{\psi^* h}
$$
on $M - (B \cup Sing)$, where $Sing$ is the set of singularities of singular
fibers of $\psi$, $\nabla_{\psi}$ is the symplectic connection defined by
$\psi$,  $X_h$ and $X_{\psi^* h}$ are Hamiltonian vector fields defined by
$\om_Y$ and $\om_M$ respectively, and $\tau$ is a positive function.

In other words, if there is a toric action of $T^k$ on $M$ we are looking for
the leaves of this action,  a family of (generically smooth) symplectic toric
manifolds with the same moment maps, parameterized by a symplectic toric
manifold. Note that the compatibility condition (2) is trivial  if both $(M,
\om_M), (Y, \om_Y)$ are Kahler, $\psi$ is a complex map and the complex
dimension of $Y$ is one.

As it is proven in \cite{4}, if a pseudotoric structure $(f_1, ..., f_k, B,
\psi, (Y, \om_Y))$ is fixed on a symplectic manifold $(M, \om_M)$ then the
choice of  moment maps $(h_1, ..., h_{n-k})$ on $Y$ defines a lagrangian
fibration on $M$.

Now we are ready to prove Proposition 1.1. Consider the direct product $\C
\PP^2_x \times \C \PP^2_y$, fix homogenious coordinates $[x_0: x_1: x_2]$ and
$[y_0: y_1: y_2]$ on the projective planes and realize $F^3$ as a hypersurface
in $\C \PP^2_x \times \C \PP^2_y$ defined by the equation
$$
\sum_{i=0}^2 x_i y_i = 0.
$$
Consider the following map:
$$
\psi: \C \PP^2_x \times \C \PP^2_x - B \to \C \PP^2_w,
$$
defined by the equations
$$
w_i = x_i y_i, i = 0, ..., 2,
$$
where $[w_0: w_1 : w_2]$ are homogenious coordinates on the last $\C \PP^2_w$.
Since the coordinates are fixed, a standard Fubini - Study metric is fixed as
well, and we endow the last $\C \PP^2_w$ by the corresponding symplectic form.
The base set $B$ consists of 6 lines defined  either by
$$
x_i = x_j = y_k = 0
$$
or by
$$
x_i = y_j = y_k = 0,
$$
where $i, j, k$ are distinct numbers from $\{0, 1, 2 \}$.

What are the fibers of $\psi$? It's not hard to see that for each $p \in \C
\PP^2_w$ with coordinates $[w_0: w_1 : w_2]$ such that $w_i \neq 0, i =0, 1, 2$
the compactified fiber
$$
\overline{\psi^{-1} (p)} = \psi^{-1}(p) \cup B
$$
is a del Pezzo surface isomorphic to $\C \PP^2_3$. The base set $B$ is a
complex six - angle on $\C \PP^2_3$. This fiber is generic.

If now we take $p \in \C \PP^2_w$ with coordinates $[w_0: w_1: w_2]$ such that
exactly one $w_i$ equals to zero then the corresponding compactified fiber
$\overline{\psi^{-1}(p)}$ is isomorphic to the union of two del Pezzo surfaces
of the form $\C \PP^2_1$ which intersect each other at a "diagonal" of the six
- angle $B$. Each $\C \PP^2_1$ has basic 4 - angle given by three edges of $B$
and the intersection line. Thus the corresponding compactified fiber is
singular with the singular set given by the "diagonal" line.

For a point $p \in \C \PP^2_w$ with coordinates $[w_0: w_1: w_2]$ such that
only one $w_i$ is non zero the compactified fiber is again non smoth --- it
consists of 2 $\C \PP^2$ and two quadrics. The six - angle $B$ is divided into
three pieces --- two conjugated triangles and two lines which form a 4 - angles
being taken with the diagonal lines. Triangles lie on $\C \PP^2$s and 4 -
angles on the quadrics.

Topologically each fiber is toric, and the construction endows each fiber with
a divisor from the anticanonical system which can play the role of the
degeneration simplex of the corresponding moment maps (or boundary divisor).
Let us show that $\C \PP^2_x \times \C \PP^2_y$ carries a pair of moment maps
whose Hamiltonian action preserves each fiber, so their restriction to each
fiber gives a completely integrable system.

For each $\C \PP^2_{*}$ in the direct product we can consider a real function
of the form
$$
f (x) = \frac{\sum_{i=0}^2 \la^{x}_i \vert x_i \vert^2}{\sum_{i=0}^2 \vert x_i
\vert^2},
$$
where $\la_i^{x}$ are distinct real numbers. In the same form one constructs a
function $f(y)$ substituting $y$ instead of $x$ in the formula above. Then
$f(x), f(y)$ are real Morse functions on the corresponding $\C \PP^2_{*}$,
whose Hamiltonian vector field preserves the Kahler structure. For two natural
projections
$$
q_*: \C \PP^2_x \times \C \PP^2_y \to \C \PP^2_*,
$$
where $*$ is either $x$ or $y$, a function of the form
$$
F = q^*_x f(x) + q^*_y  f(y)
$$
(where now the star means the lifting) is a real Morse function on the direct
product, whose Hamiltonian vector field preserves the Kahler structure on $\C
\PP^2_x \times \C \PP^2_y$, and it is not hard to see that if the following
condition on $\la_i^x, \la_j^y$ holds
$$
\la^x_0 + \la_0^y = \la^x_1 + \la^y_1 = \la^x_2 + \la^y_2
$$
then the Hamiltonian vector field of such a function $F$ must preserve the
fibers of $\psi$. Since we have on $\C \PP^2$ exactly two algerbaically
independent  almost everywhere real Morse  functions of this form, it follows
that there are exactly two algebraically independent almost everywhere Morse
functions $F_1, F_2$ which preserve the fibers. The degeneration simplex
$$
\De(F_1, F_2) = \{ X_{F_1} \wedge X_{F_2} = 0 \} \subset \C \PP^2_x \times \C
\PP^2_y
$$
consists of the six - angle $B$ and three its "diagonal" lines. Thus according
to the general theory, see \cite{4}, the degeneration simplex is exactly the
union of $B$ and $\Sing$.

Now we claim, that the constructed above data $(F_1, F_2, B, \psi, \C \PP^2)$
is a {\it pseudotoric structure} on the direct product $\C \PP^2_x \times \C
\PP^2_y$. To prove it remains to show that the second compatibility condition
from the definition holds, which is not hard to do and we leave it aside the
discussion. Our aim is to show that $F^3$ is pseudotoric as well.

The desired pseudotoric structure is given by the restriction of the data
$(F_1, F_2, B, \psi, \C \PP^2)$ to the line $\C \PP^1_w \subset \C \PP^2_w$
with the equation
$$
w_0 + w_1 + w_2 = 0.
$$
Then since for $F_1, F_2$ their Hamiltonian vector fields are tangent to $F^3$
and the image of $F^3$ under the map $\psi$ is exactly the line $\C \PP^1_w$ we
have the data $(f_1, f_2, B, \psi_1, \C \PP^1)$ where

$f_i = F_i|_{F^3}$;

$B = B$ since $B$ lies on $F^3$;

$\psi_1 = \psi|_{F^3 - B}$;

and it's not hard to see that the data $(f_1, f_2, B, \psi_1, \C \PP^1_w)$ is a
pseudotoric structure on $F^3$.

The line $ \C \PP^1_w = \{ w_0 + w_1 + w_2 = 0 \} \subset \C \PP^2$ intersects
three basic lines $l_i = \{ w_i = 0 \}, i =0, 1, 2,$ at three points, and these
points underlay singular fibers of $\psi$. The type of the singular fibers was
described above: each singular fiber is a pair of two $\C \PP^2_1$ (del Pezzo
surfaces of degree 8) with the singular set represented by a "diagonal" line.

Now we are ready to present the minimal lagrangian fibration of $F^3$,
mentioned in the Introduction. To get it let us fix a smooth Morse function on
$\C \PP^1$ which has two critical points, maximal and minimal, say, at $[0: 1:
-1]$ and $[1: 0: -1]$. Then the boundary divisor $D \subset F^3$ of the
corresponding fibration consists of two singular fibers, and each of these
fiber is formed by two del Pezzo surfaces $\C \PP^2_1$. Thus totally it gives 4
del Pezzo surfaces which are toric with the intergals $F_1, F_2$, restricted to
each component, and these integrals induce a canonical toric fibration of each
component. It shows that $D$ is fibered by smooth tori as usual in toric
geometry. Now consider what happens on $F^3 - D$. Each fiber $\psi^{-1}(p)$ is
fibered on smooth two dimensional tori by the restriction of $f_1$ and $f_2$.
These tori are lebeled by two regular values of $f_1$ and $f_2$. Take a smooth
level loop $\ga$ of the function $h$ on $\C \PP^1$ which doesn't pass through
the point $[1: -1: 0]$. Fix a pair of regular values for $f_1, f_2$ and collect
the corresponding tori with the same values along $\ga$:
$$
\cup_{p \in \ga} T^{c_1, c_2}_p,
$$
where $T^{c_1, c_2}_p$ is the torus from the fiber over $p$ with the values
$c_1, c_2$. Then this union is a smooth three dimensional torus which is
lagrangian in $F^3$. This is a generic fiber.

Now when singular tori appear? If $\ga$ passes through the singular point $[1:
-1: 0]$ it doesn't mean automatically that the corresponding torus must be
singular. Indeed for the inner points of the six - angle in $\R^2$ which is the
set of values for the pair $f_1, f_2$ a torus over $[1: -1: 0]$ is either 2 -
dimensional or 1 - dimensional if the values $(c_1, c_2)$ lie on the diagonal
segment in the six - angle. Generic choice of $(c_1, c_2)$ gives a 2 -
dimensional torus, and collecting tori along $\ga$ we get a smooth lagrangian
tori for this general case. The singularity appears if the pair $(c_1, c_2)$
lies on the diagonal segment, and then this segment parameterizes singular tori
of the following type: $C^2_p$ is applied to $T^3$.

The lagrangian fibration of $F^3$ presented above is minimal in the following
sense: another choice of a function $h$ on $\C \PP^1$ underlying $\psi_1$ gives
another fibration with bigger set of singular tori (except this $h$ satisfies
the same properties so its maximal and minimal points coincide with any two of
the three marked points).

\section{Specialty condition}

Recall, that first time the notion of special lagrangian tori was presented for
Calabi -- Yau manifolds as a property, "dual" to the holomorphity condition in
the framework of the Mirror Symmetry problem. D. Auroux generalizes it in
\cite{2} for the case of certain Fano varieties. Namely, see \cite{2}, if a
Fano variety $X$ admits an effective divisor in the anticanonical system, say,
$D \subset \vert - K_X \vert$, then there exists a top holomorphic form
$\theta_D$ on the complement $X - D$; this form has pole at $D$ and doesn't
vanish on $X -D$. Note that this form is uniquelly defined by $D$ up to $\C^*$.

Since by the very definition a Fano variety $X$ can be embedded to a projective
space by the anticanonical system one can attach to $X$ the corresponding
symplectic form given by the restriction to the image of $X$ of the standard
symplectic form on the projective space. Thus any Fano variety carries a
symplectic structure with the symplectic form, proportional to the (anti)
canonical class, therefore any Fano variety can be treated as a monotone
symplectic manifold; the symplectic structure on $X$ can be call {\it
canonical} (or anticanonical, if one wishes).

A lagrangian torus $S \subset X$ is {\it special} in $(X, D)$, if the imaginary
part $\Im e^{i s} \theta_D|_S$ of  of the restriction $\theta_D|_S$ vanishes
for certain twist $e^{i s} \in U(1)$. A lagrangian fibration of $X - D$ is
called {\it special} if the condition holds for any fiber with the same twist
$e^{i s}$.

The simplest example of special lagrangian fibration comes from the toric
geometry consideration. Indeed, the canonical lagrangian fibration of a toric
manifold $X$ degenerates dimensionally at a divisor $D_0$ from the
anticanonical system $\vert - K_X \vert$; then  the corresponding holomorphic
top form $\theta_{D_0}$ has very nice form in the action - angle coordinates
which implies the property.

But the standard fibration is not quite good for several reasons, f.e. the
boundary divisor $D_0$ is always reducible, and it was the problem stated in
[2] to find a special lagrangian fibration for a smooth element from the
anticanonical system on $\C \PP^2$. The example of non standard special
lagrangian fibration given in \cite{2} has a reducible boundary divisor, the
sum of line and conic, but it is more generic than the toric one.

At the same time the problem can be studied over any Fano variety with an
effective divisor in the anticanonical system. F.e., we claim that the flag
variety $F^3$ admits a special lagrangian fibration, and this fibration is our
minimal one, constructed above.

\begin{prop} The minimal lagrangian fibration on $F^3$ is special lagrangian.
\end{prop}

Indeed, the boundary divisor of the minimal fibration consists of two
compactified fibers of $\psi$. The anticanonical class of $F^3$ as a bundle is
given by the product $q_x^* \Oh_{\C \PP^2_x}(2) \otimes q_y^* \Oh_{\C \PP^2_y}
(2)$ by the adjuction formula, and at the same time the compactified fibers of
$\psi$ are sections of the bundle $q_x^* \Oh_{\C \PP^2_x}(1) \otimes q_y^*
\Oh_{\C \PP^2_y} (1)$ and thus the boundary divisor of the minimal fibration
lies in the anticanonical system --- the homological condition is satisfied.

Take now the boundary divisor $D \subset F^3$ of the minimal fibration and
consider the corresponding top holomoprhic form $\theta_D$ on $F^3 - D$.

We claim that there exist three real vector fields on $F^3 - D$, two global and
one local, which preserves the holomorphic form $\theta_D$. This means that the
Lie derivative of $\theta_D$ along any one from this triple is zero, and since
$\theta_D$ is closed and the vector fields commute one gets a constant function
substituting these fields to $\theta_D$.

These three vector fields are: the Hamiltonian vector fields $X_{f_1}, X_{f_2}$
(defined globally on whole $F^3$ and the lifted vector field $\nabla_{\psi}
X_h$, which is defined outside of the singular points of $\psi$. Indeed, each
$X_{f_i}$ generates Hamiltonian isotopy of $F^3$ which leaves the boundary
divisor $D$ unmoved, and this implies that the Lie derivative $\sL_{X_{f_i}}
\theta_D$ is zero. On the other hand, the Hamiltonian vector field $X_{h}$ on
the base of the pseudotoric structure generates certain 1 - parameter family of
symplectomorphisms of $F^3$ which moves fibers of $\psi$ to fibers and which
leaves unmoved the boundary divisor. This means that the Lie derivative
$\sL_{\nabla_{\psi}(X_h)} \theta_D$ vanishes (but now it is defined outside of
the singular points of $\psi$ only).

Consider now any torus $S \subset F^3$ from the minimal lagrangian fibration.
Denote as $S_0 \subset S$ the "smooth" part consists of smooth points (so
generically $S_0 = S$). At each point $p \in S_0$ the tangent space $T_pS$ is
spanned by $X_{f_1}, X_{f_2}, \nabla_{\psi} X_h$. At the same time, since
$\theta_D(X_{f_1}, X_{f_2}, \nabla_{\psi} X_h)$ is a constant complex function
on $S_0$ and the determinant space $ det T_pS$ is real one dimensional, the
restriction of $\theta_D$ to $S_0$ must be proportional to any real non
vanishing top form with the proportionality coefficient of the form $f e^{i
\phi}$ where $f$ is a real function and $\phi$ is a constant. And the point is
that this constant (so the argument) is the same for all $S \subset F^3$ from
the minimal fibration. This ends the proof.

Here we implicitly use the fact that the function $h$ is a symbol on the based
$\C \PP^1_w$, so its Hamiltonian vector field preserves the complex structure.
Therefore the statement of Proposition 2.1 remains true for the  lagrangian
fibration on $F^3$ given by {\it any} symbol on the base $\C \PP^1_w$.

\section{Toric degeneration of $F^3$}

We could describe the pseudotoric structure $(f_1, f_2, B, \psi_1, \C \PP^1)$
on the flag variety $F^3$ without any references to the representation by a
hypersurface in the direct product $\C \PP^2_x \times \C \PP^2_y$, but the
computations above are usefull for the studying of toric degeneration of $F^3$
and a relationship between its canonical toric fibration and the minimal
fibration of $F^3$ constructed above.

Toric degeneration of $F^3$ is a deformation of the flag variety to a singular
toric variety. In the coordinates $\{ x_i, y_i \}$ it is given by the following
family
$$
F_t = \{ t x_0 y_0 + x_1 y_1 + x_2 y_2 = 0 \} \subset \C \PP^2_x \times \C
\PP^2_y,
$$
where $t$ is a complex parameter. When $t \neq 0$ the variety $F_t$ is
isomorphic to $F^3$ while $t=0$ implies that $F_0$ is a singular toric variety.
As a toric variety $F_0$ carries a canonical lagrangian fibration. Our aim for
this section is to relate the fibers of the lagrangian fibration of $F^3$
constructed in the previous section and the fibers of this canonical lagrangian
fibration.

The first observation which is usefull in this way is the following
\begin{prop} The singular toric variety $F_0$ carries a pseudotoric structure, and
the standard toric fibration can be realized by the choice of an appropriate
function $h$ on the base $\C \PP^1$.
\end{prop}

Indeed, consider again $\C \PP^2_x \times \C \PP^2_y$ and the pseudotoric
structure $(F_1, F_2, B, \psi, \C \PP^2_w)$, presented in a previous section.
Now we take the line
$$
\C \PP^1_0 = \{w_1 + w_2 = 0 \} \subset \C \PP^2_w;
$$
it is clear that the Hamiltonian actions of $F_1$ and $F_2$ both preserve
$F_0$, therefore the restriction of $F_i$ to $F_0$ gives a set of commuting
functions, $B$ is the base set of the family, the restriction of $\psi$ to
$F_0$ gives a map to $\C \PP^1_0$ which we denote as $\psi_0$, and totally it
gives a pseudotoric structure on $F_0$.

The singular fibers of $\psi_0$ are the fibers over two points: $[1: 0: 0]$ and
$[0: 1: -1]$ in $\C \PP^1_0$. Now one can restore a standard toric fibration of
$F_0$:

\begin{prop} Taking function $h$ on $\C \PP^1_0$ such that it has exactly two critical
points,  at $[1: 0: 0]$ and $[0: 1: -1]$ one gets a canonical toric fibration
of $F_0$.
\end{prop}

 One can show it in two steps:

--- first, it's clear that a standard toric fibration on $F_0$ can be constructed from
a standard toric fibration on the ambient variety $\C \PP^2_x \times \C
\PP^2_y$ which is invariant with respect to $F_1$ and $F_2$;

--- second, one can show that a standard toric fibration of the product variety
$\C \PP^2_x \times \C \PP^2_y$ is defined by certain choice of moment maps
$H_1, H_2$ on $\C \PP^2_w$ and the pseudotoric structure $(F_1, F_2, B_0, \psi,
\C \PP^2_w)$.

These $H_1, H_2$ must be diagonal in the coordinates $[w_0: w_1: w_2]$ moreover
they can be choosen such that one of them, say, $H_1$ is constant on the line
$l_0 = \{ w_0 = 0 \}$. Then the Hamiltonian vector field of $H_1$ must be
parallel to the line $\C \PP^1_0$, and the restriction of $H_1$ to $\C \PP^1_0$
would be our $h$.

The only essential part of the proof is to show that diagonal in the
coordinates $[w_0: w_1 : w_2]$ functions $H_1$ and $H_2$ give a standard toric
fibration of $\C \PP^2_x \times \C \PP^2_y$, so we clarify here only this
point. To see this note that for the choice $H_1, H_2$ diagonal the singular
fibers of $\psi$ lie over the degeneration simplex $\De(H_1, H_2)$ of functions
$H_1, H_2$ which is the union of lines $l_i = \{ w_i = 0 \} \subset \C
\PP^2_w$. This means that the boundary divisor $D$ of the corresponding
fibration on $\C \PP^2_x \times \C \PP^2_y$ consists of the components:
$$
\overline{\psi^{-1}(l_i )}, i = 0, 1, 2,
$$
and these are presicely six irreducible components of the form $q_x^*(l^x_i),
q_y^* l^y_i, i = 0, 1, 2,$ where $l_i^* \subset \C \PP^2_{*}$ is the projective
line
$$
l^x_i = \{ x_i = 0 \}, l^y_i = \{ y_i = 0 \}, i = 0, 1, 2.
$$
Thus the boundary divisor for the fibration, generated by $H_1, H_2$, is the
same as for a standard toric fibration, and since for the diagonal functions
$H_1, H_2$ the corresponding fibration has smooth fibers only, this fibration
must be isomorphic to a standard fibration (given by the product of two
Clifford fibrations on both $\C \PP^2_x, \C \PP^2_y$.

How one can compare lagrangian fibrations on $F_0$ and $F^3 = F_1$? The regular
fibers can be related by certain Hamiltonian trasnformations of $\C \PP^2_x
\times \C \PP^2_y$. The procedure can be done as follows. Consider a Morse
function $g$ on $\C \PP^2_w$ such that its Hamiltonian vector field moves in a
finite time the projective line $l_1 = \{w_0+  w_1 + w_2 = 0 \}$ to the
projective line $l_0 = \{w_1 +w_2 \}$. If a Morse function $h_1$ is fixed on
the  projective line $l_1$ then it can be deformed to a Morse function $h_0$ on
$l_0$. At the same time one can define a Hamiltonian isotopy on $\C \PP^2_x
\times \C \PP^2_y$ generated by a partial lift of $g$. Since $\psi^* g$ is not
defined at the base set $B \subset \C \PP^2_x \times \C \PP^2_y$ and $\psi^* g$
is not smooth at the singular set $\Sing$ we take two small neighborhoods $B_1
\supset B_2$ of
$$
\De(F_1, F_2) = B \cup \Sing \subset \C \PP^2_x \times \C \PP^2_y
$$
of sufficiently small radii $1 >> r_1 > r_2 > 0$; then one can construct a
smooth function $G$ on whole $\C \PP^2_x \times \C \PP^2_y$ by the following
rules:

--- $G = \psi^* g$ on $\C \PP^2_x \times \C \PP^2_y - B_1$;

--- $G$ is zero in $B_2$;

--- $G$ smoothly changes from $\partial B_1$ to $\partial B_2$.

Then it is clear that the Hamiltonian action of $G$ on $\C \PP^2_x \times \C
\PP^2_y$ moves a ``central part'' $F_1^c$ of $F_1$ to a ``central part''
$F^c_0$ of $F_0$ where one means the ``central parts'' to be
$$
F^c_i = F_i - (F_i \cap B_1), i =0, 1.
$$
in a finite time. Therefore one can state the following

\begin{prop} Every smooth lagrangian torus of the minimal lagrangian fibration of
$F^3$ constructed in the previous section is Hamiltonically isotopical to a
smooth lagrangian torus of the standard lagrangian fibration of the toric
degeneration $F_0$.
\end{prop}

To prove this fact one takes $r_1$ such that a given torus doesn't intersect
$B_1$; it's easy  to see that an appropriate choice exists. Then the
Hamiltonian isotopy generated by $G$ moves our given torus to a submanifold of
$F_0$ which is fibered over the corresponding level loop of $h_0$; by the
construction of $G$ this function commutes with $F_i$ everywhere outside of
$B_1$ which implies that the resulting torus in $F_0$ lies on a mutual level
set of the functions, and it follows that it is lagrangian in $F_0$. It remains
to see that it is given by the pseudotoric construction for the function $h_0$
which has exactly two critical points underlying singular fibers of $\psi_0$,
and hence it is a fiber of the standard fibration.

A singular torus via the limiting procedure $t \to 0$ collapses to an
isotropical torus in $F_0$ which lies on a boundary component.

The last proposition implies the following fact:

\begin{cor} The Floer cohomology $FH(S)$ of a smooth  fiber $S \subset F^3$ in the minimal
lagrangian fibration is the same as the Floer cohomology of a fiber of the
standard lagrangian fibration of $F_0$, the toric degeneration of $F^3$.
\end{cor}

Indeed, the Floer cohomology (whatever it is defined) must be stable with
respect to the Hamiltonian isotopy, and all the ingredients exploited for
calculations of $FH(S)$ on $F_0$ can be moved to $F^3$ giving the same answer.

\section{Digression: geometric construction}

As we've mentioned above the pseudo toric structure on $F^3$ can be described
in pure geometric terms. Here we discuss briefly how it can be done since
probably it will be usefull for generic flag variety.

Consider $\C^3$ and fix a basis $e_1, e_2 , e_3$ and a hermitian product $<,>$
such that $<e_i>$ is an orthonormal basis for $<,>$. Consider the
projectivization $\C \PP^2$ of $\C^3$; the hermitian product defines a standard
Kahler triple on $\C \PP^2$, and the basis $<e_i>$ fixes a "basic" triangle
$\De \subset \C \PP^2$. Consider the flag variety $F^3$; the hermitian product
$<,>$ defines a Kahler triple on $F^3$ and the basis $<e_i>$ gives six points
on $F^3$ which correspond to ordered pairs $(e_i, e_j), i \neq j$; there are
six edges which relate these six points such that one has a six - angle $\hat
\De \subset F^3$.

Any self adjoint operator $A$ on $(\C^3, <,>)$ induces a pair of real functions
defined up to constant on $\C \PP^2$ and $F^3$. Namely $A$ generates an
infinitesimal symmetry of $(\C^3, <,>)$ and thus the corresponding
infinitesimal symmetries of $\C \PP^2$ and $F^3$ preserve the Kahler structures
and therefore they are Hamiltonian and there are real functions $f_A, F_A$ on
$\C \PP^2$ and $F^3$ respectively defined up to constant. If one takes the
standard projection
$$
\pi: F^3 \to \C \PP^2
$$
which is the forgetfull map for the second element in the flag, then it's not
hard to see that $\pi$ is a symplectic map so for each pair of tangent vectors
$v_1, v_2 \in T_p F^3$ the value of the symplectic form $\om_{F^3}(v_1, v_2)$
equals to $\om_{\C \PP^2}(\pi(v_1), \pi(v_2))$ up to a multiple constant. It's
not hard to see that the functions $f_A$ and $F_A$ are related by $\pi$ as
follows:
$$
X_{f_A} = \pi(X_{F_A}).
$$

Consider a pair of self adjoint operators $A_1, A_2$ which are diagonal with
respect to the basis $<e_i>$ with distinct eigenvalues, such that they are not
proportional and their sum is not proprotional to identity. Then they give two
pair of Morse functions: $f_1, f_2$ on $\C \PP^2$ and $F_1, F_2$ on $F^3$ such
that
$$
X_{f_i} = d \pi(X_{F_i}).
$$
The degeneration simplices
$$
\De(f_1, f_2) = \{ X_{f_1} \wedge X_{f_2} = 0 \} \subset \C \PP^2
$$
and
$$
\De(F_1, F_2) = \{ X_{F_1} \wedge X_{F_2} = 0 \} \subset F^3
$$
are related in a similiar manner
$$
\pi(\De(F_1, F^3)) = \De(f_1, f_2),
$$
and the first degeneration simplex $\De(F_1, F_2)$ consists of such flags $(p,
l)| p \in l$ that the Hamiltonian actions of $X_{f_1}$ and $X_{f_2}$ on $(p,
l)$ are proportional. Therefore it's clear that $\hat \De$ is a piece of $\De
(F_1, F_2)$; it remains to add three "diagonal" lines in the six - angle $\hat
\De$ to get $\De(F_1, F_2)$.

The choice of $A_1, A_2$ defines a singular connection $\nabla_A$ on the bundle
$\pi: F^3 \to \C \PP^2$; the horizontal distribution over a point $p \in \C
\PP^2$ which doesn't belong to $\De$ is given by the linear span of $X_{F_i},
I(X_{F_i})$ where $I$ is the complex structure. This connection degenerates at
$\pi^{-1}(\De)$. And a natural question arises --- does this singular
connection admits covariantly constant sections? Or is this distribution
integrable? The integrability means the following: take a fiber of $\pi$ and
for each point of the fiber there is a (generically smooth) complex submanifold
of $F^3$ which is

--- projected to the base $\C \PP^2$ by the map $\pi$;

--- over each point except $\De$ the tangent space to this submanifold is horizontal
with respect to the singular connection $\nabla_A$.

Thus the section can meet each other over the degeneration set $\De$ and not
outside.

Now we show that it is indeed the case:

\begin{prop} The distribution, defined by the singular connection $\nabla_A$, is intergable.
Thus the connection $\nabla_A$ admits covariantly constant sections, and the
space of the sections is parameterized by the fiber $\pi^{-1}(p)$ for a generic
point $p \in \C \PP^2$.
\end{prop}

The proof is strightforward. We fix a point $p \in \C \PP^2$ outside of $\De$
and consider the fiber of $\pi$ over $p$. The fiber is the pencil of lines
passing through $p$. If we fix an  element from the pencil, a line $l \subset
\C \PP^2$, then there are two possibilities:

--- either $l$ doesn't intersect vertices of $\De$, and this case is generic;

--- or $l$ passes through a vertex of $\De$, and there are precisely three such elements in the pencil.

First we fix a generic line $l \ni p$ from the pencil and construct the
corresponding covariantly constant section. For a fixed pair $(p, l) | p \in
l,$ consider the action of all non degenerated complex operators, diagonal in
the basis $<e_i>$, up to scale, so the action of a complex 2 -torus. The orbit
is isomorphic to $\C \PP^2 - \De$; denote it as $\psi^{-1}(p)$. Consider the
degeneration simplex $\De(F_1, F_2)$ for the functions $F_1, F_2$ in $F^3$. It
consists of

(i) 6 vertices which correspond to pairs {\it (vertex of $\De$, edge of
$\De$)},

(ii) six edges of two types, {\ it (vertex of $\De$, a line through this
vertex)} and {\it (point on edge of $\De$, the edge of $\De$)},

(iii) and three diagonal lines {\it (point on edge of $\De$, line through
vertex of $\De$)} such that the point lies on the line.

The degeneration simplex $\De(F_1, F_2)$ is decomposed on two parts: the "six -
angle" $B$ is formed by (i) and (ii) and the subset $\Sing$ is (iii). The
division is defined by the following argument: the six - angle $B$ is mapped to
a convex polytop in $\R^2$ by the action map
$$
(F_1, F_2): F^3 \to \R^2,
$$
while the image of $\Sing$ under this map lies inside of the convex polytop.
Now it is not hard to see that the union
$$
\psi^{-1}(p) \cup B
$$
is a smooth submanifold of $F^3$ birationally isomorphic to $\C \PP^2$ (since
the projection $\pi$ is a birational map). Moreover, the union
$$
\psi^{-1}(p) \cup B
$$
is isomorphic to  del Pezzo surface $\C \PP^2_3$. It is toric; moment maps are
gien by the restriction of $F_1, F_2$ to it. If one wants it can be established
more carefully, step by step.

 {\it The first step: operators of rank 3.} It's clear that the
projective action of the diagonal operators of rank 3 on the flag $(p, l)$ is
transitive; each point of the base $\C \PP^2$ outside of $\De$ can be reached,
and the second element of the pair $(K(p), K(l))$ is uniquelly defined by the
first one, where $K$ is an operator. Thus one gets the part of orbit isomorphic
to $\C \PP^2 - \De$.

{\it The second step: operators of rank 2.} If $K$ has one zero eigenvalue then
$(K(p), K(l))$ is of the following type: $K(p)$ lies on the corresponding edge
of $\De$ and $K(l)$ is this edge. At the same the action of $K(p)$ is
transitive, and this part of orbit is isomorphic to the union of three lines
without the intersection points.

{\it The third step: operators of rank 1.} At this step the picture degenerates
since $l$ vanishes under the action of such an operator. At the same time point
$p$ goes to a vertex of $\De$, so totally our orbit is projected to whole $\C
\PP^2$. But upstairs in $F^3$ the degeneration leads to the uncertainty at
these vertices: the limiting procedure along  edges gives different limiting
lines at the vertex where the edges meet each other. This means that the
appropriate smooth compactifiction must be given by the blow up of $\C \PP^2$
at the vertices of $\De$.  And this is our covariantly constant section which
is isomorphic to the del Pezzo surface $\C \PP^2_3$ embedded to $F^3$
horizontally. By the construction it is parallel to the Hamiltonian vector
fields $X_{F_i}$ almost everywhere.

It's clear that the horizontal section we've constructed starting with a
general flag $(p, l)$ contains the six - angle $B \subset F^3$  formed by edges
and vertices of $\De \subset \C \PP^2$.

The section endowed with the restricted moment maps $F_i$ gives a symplectic
toric manifold, and the six - angle is precisely the boundary divisor. Note
that each general horizontal section contains this basic six - angle by the
construction.

Now consider the special case when the line $l$ from the pair $(p, l)$ passes
through a vertex of $\De$, say, the point $[1:0:0]$ in the homogenious
coordinates $[z_0: z_1: z_2]$ on the base $\C \PP^2$. Take again diagonal
operators and study their action on $(p, l)$. Let's again study this case step
by step and then put down certain generic arguments.

{\i The first step: operators of rank 3.} Again by the action of any such
operator we can get any point outside of $\De$ in $\C \PP^2$. At the same time
all the lines we get will pass through the point $[1:0:0]$. It shows that this
part is isomorphic to $\C \PP^2 - \De$.

{\it The second step: operators of rank 2.} Here we reach certain difference
with the generic case since if $K$ has the first eigenvalue equals to zero then
the pair $(K(p), K(l))$ degenerates. Indeed, the line $l$ should be transformed
to the line $\{z_0 = 0 \}$ but on the other hand it should pass through the
point $[1:0:0]$ --- and thus it vanishes. To complete this part one can restore
the second element of the pair taking the line passing through $[1:0:0]$. Note
that the line of $F^3$ which corresponds to the last case is a part of the
degeneration simplex of $f_1, f_2, f_3$ on $F^3$.

{\it The third step: operators of rank 1.} Near the point $[1:0:0]$ we have the
same picture as in the general case so it must by blown - up. But now the other
vertices are smoothly lifted to $F^3$. Therefore this part of the compactified
section is isomorphic to $\C \PP^2_1$ blown up at the point $[1:0:0]$.

 At the same time one meets certain discrepancy: the constructed smooth section doesn't
contain the basic six - angle $B$. Indeed, by the construction our section
doesn't contain two vertices of $B$. At the same time it contains a diagonal
line of $B$ where the system of integrals $f_1, f_2, f_3$ degenerates.  To
recover the property that each element of the family contains a base set one
needs

{\it The fourth step: adding a component.} Note that the component we get above
is a Schubert submanifold of $F^3$. Indeed, it is presented by the condition
that line $l$ from  flag $(p, l)$ passes through a fixed point. Denote it by
$D_{p_0} \subset F^3$, where $p_0$ is the point $[1:0:0]$. The diagonal line of
the six - angle $B$ at the same time is a Schubert submanifold as well: it
consists of the flags $(p, l)$ such that $p$ lies on the line $\{z_0 = 0 \}$
and $l$ passes through the point $[1:0:0]$. The additional component one needs
to take is "dual" to $D_{p_0}$: it is a Schubert submanifold consists of the
flags $(p, l)$ such that $p$ lies on the line $ l_0 = \{ z_0 = 0 \}$. We denote
it as $D_{l_0} \subset F^3$. Now it's clear that $D_{x_0} \cup D_{l_0}$ is a
singular submanifold, invariant under the Hamiltonian action of $f_1, f_2,
f_3$, which contains the base six - angle $B$. And it is our singular
horizontal section, which corresponds to the element of the fiber $\pi^{-1}(p)$
passing through $p_0$.

Since we have described all the possibilities for  lines through $p$ this ends
the proof.

\end{document}